\documentclass[9pt,technote]{IEEEtran}
\usepackage{amsthm}
\usepackage{amsmath}
\usepackage{amssymb}
\usepackage{amsfonts}
\usepackage{graphicx}

\newtheorem{assumption}{Assumption}

\makeatletter
\theoremstyle{plain}
\newtheorem{thm}{\protect\theoremname}
\theoremstyle{remark}
\newtheorem{rem}{\protect\remarkname}
\theoremstyle{definition}
\newtheorem{problem}{\protect\problemname}

\theoremstyle{definition}
\newtheorem{example}{\protect\examplename}
\theoremstyle{plain}
\newtheorem{lem}{\protect\lemmaname}
\theoremstyle{plain}
\newtheorem{cor}{\protect\corollaryname}
\theoremstyle{plain}
\newtheorem{prop}{\protect\propositionname}
\providecommand{\corollaryname}{Corollary}
\providecommand{\lemmaname}{Lemma}
\providecommand{\problemname}{Problem}
\providecommand{\examplename}{Example}
\providecommand{\remarkname}{Remark}
\providecommand{\theoremname}{Theorem}
\providecommand{\propositionname}{Proposition}

\IEEEoverridecommandlockouts                             

\newcommand{\R}{{\mathbb R}}

\newcommand{\M}{{\mathbb M}}
\newcommand{\Z}{{\mathbb Z}}

\newcommand{\tr}{{\text{tr}}}

\begin{document}
\date{}
\title{\Large On the Optimal Solutions of the Infinite-Horizon \\Linear Sensor Scheduling Problem}
\author{\normalsize Lin Zhao, Wei Zhang, Jianghai Hu, Alessandro Abate and Claire J. Tomlin
\thanks{L. Zhao and W. Zhang are with the Department of Electrical and Computer Engineering, The Ohio State University (e-mail:~\{zhaol, zhang\}@ece.osu.edu)}
\thanks{J. Hu is with the School of Electrical and Computer Engineering, Purdue University, West Lafayette, IN (e-mail:~jianghai@purdue.edu)}
\thanks{A. Abate is with the Department of Computer Sciences, University of Oxford, United Kingdom, and
with the Delft Center for Systems \& Control, Delft University of Technology, The Netherlands (e-mail:alessandro.abate@cs.ox.ac.uk)}
\thanks{C. J. Tomlin is with the Department of Electrical Engineering and Computer Sciences, University of California at Berkeley (e-mail:~tomlin@eecs.berkeley.edu)}
}
 
\maketitle
\begin{abstract}
This paper studies the infinite-horizon sensor scheduling problem
for linear Gaussian processes with linear measurement functions. Several
important properties of the optimal infinite-horizon schedules are
derived. In particular, it is proved that under some mild conditions,
both the optimal infinite-horizon average-per-stage cost and the corresponding
optimal sensor schedules are independent of the covariance matrix
of the initial state. It is also proved that the optimal estimation
cost can be approximated arbitrarily closely by a periodic schedule
with a finite period. Moreover, it is shown that the sequence of the 
average-per-stage costs of the optimal schedule must
converge. These theoretical results provide valuable insights into the design and
analysis of various infinite-horizon sensor scheduling algorithms.
\end{abstract}

\section{Introduction}

The sensor scheduling problem seeks an optimal schedule over a certain
time horizon to activate/deactivate a subset of available sensors
to improve the estimation performance and reduce the estimation cost
(e.g. energy consumption and communication overheads). It has numerous
applications in various engineering fields~\cite{LiDSP09,MeiTAC67,MOTR06}.

Previous research has mainly focused on the finite-horizon sensor
scheduling problem for linear Gaussian processes. In this case, a
straightforward solution is to enumerate all the possible finite-horizon
schedules~\cite{MeiTAC67}. The complexity of such an approach grows
exponentially fast as the horizon size increases. Various methods
have been proposed in the literature to tackle this challenge. These
methods can be roughly divided into the following three categories:
(i) methods that focus on certain simple special classes of schedules,
such as myopic schedules that only consider immediate performance
at each time step instead of the overall performance over the whole
horizon~\cite{LOACC99,OSTAES}; (ii) methods that ``embed'' discrete
schedules into a larger class of schedules with continuously-variable
sensor indices~\cite{GupAuto06,NYACC09}; (iii) and methods that
prune the search tree based on certain properties of the Riccati recursions~\cite{FEAUTO08,ACC10_ss}.

The methods in the first category are often easy to implement, but
provide no guarantees on the overall estimation performance. The ``embedding''
approach in the second category is a common trick to tackle complex
discrete optimization or optimal control problems~\cite{BeAUTO05,JOTSP09}.
The resulting relaxed schedule can often be interpreted as the time-average
``frequencies'' or ``probabilities'' for using different sensors.
It has been recently proved that, in continuous time, the performance
of the optimal relaxed schedule can be approximated with arbitrary
accuracy by a discrete schedule through fast switchings~\cite{NYACC09}.
This is analogous to the result derived in~\cite{BeAUTO05} for solving
the optimal control problem of switched systems using embedding. However,
in discrete time, the result no longer holds as the switching rate
is fixed; in this case, the relaxed schedule can only be implemented
probabilistically~\cite{GupAuto06}, resulting in a random scheduling
of the sensors. The pruning methods in the third category make essential
use of the monotonicity and concavity properties of the Riccati mapping
to obtain conditions under which the exploration of certain branches
can be avoided without losing the optimal schedule. In our earlier
paper~\cite{ACC10_ss}, an efficient sub-optimal algorithm was proposed
to prune out not only the non-optimal branches but also less important
ones to further reduce the complexity. Further error bounds associated
with this pruning algorithm have also been derived in~\cite{CDC10_Mike}. In contrast to these offline scheduling methods, online event-based sensor scheduling problems have also been studied recently~\cite{Lemmon2010,Shi2011, Sijs2009}. It has been shown that these event-based online approaches can improve the performance if communication overheads are incorporated into the cost function or the constraints of the scheduling problem.

Different from most previous research, this paper studies the infinite-horizon
sensor scheduling problem for discrete-time linear Gaussian processes
observed by linear sensors. The problem is much more challenging than
its finite-horizon counterpart and has not been adequately investigated
in the literature. Instead of proposing a specific scheduling algorithm,
we focus on deriving several properties of the problem which are of
fundamental importance for the design and analysis of various infinite-horizon
sensor scheduling algorithms. In particular, it is proved that under
some mild conditions, both the optimal infinite-horizon average-per-stage
cost and the corresponding optimal sensor schedule are independent
of the covariance matrix of the initial state. It is also proved that
the optimal estimation cost can be approximated arbitrarily closely
by a periodic schedule with a finite period. Furthermore, it is concluded
that the sequence of the average-per-stage costs of the optimal schedule
must converge. These theoretical properties provide us valuable insight
into the infinite-horizon sensor scheduling problem and will be useful
for developing algorithms. In addition, the existence of a periodic
suboptimal schedule justifies the experimental results of many finite-horizon
scheduling algorithms~\cite{HoCDC07,ACC10_ss} that yield periodic
schedules for relatively large horizons.

The rest of the paper is organized as follows. The infinite-horizon
sensor scheduling problem is formulated in Section~\ref{sec:prob}.
Some important properties of the difference Riccati recursion are
reviewed in Section~\ref{sec:srm}. These properties are then used
in Section~\ref{sec:pos} to prove the universal approximation property
of the periodic schedule. Finally, some concluding remarks are given
in Section~\ref{sec:con}.

\textbf{Notation:} Let ${\cal A}$ be the semi-definite cone, namely,
the set of all the positive semidefinite matrices. Denote by $\lambda_{\min}(\cdot)$
and $\lambda_{\max}(\cdot)$ the smallest and the largest eigenvalues,
respectively, of a given matrix in ${\cal A}$. Let $\R_{+}$ and
$\Z_{+}$ be the set of nonnegative real numbers and integers, respectively.
Denote by $|\cdot|$ the standard Euclidean norm of vectors or absolute
value of numbers, and $\|\cdot\|$ the vector-induced matrix norm.
For any $\phi_{c}\in{\cal A}$ and $r>0$, define ${\cal B}(\phi_{c};r):=\{\phi\in{\cal A}:\|\phi-\phi_{c}\|\le r\}$.
Denote by $I_{n}$ the identity matrix of dimension $n$, and diag$\{.,.\}$
the diagonal matrix composed of the input arguments.

\section{Problem Formulation \label{sec:prob}}

Consider the following linear time-invariant stochastic system: 
\begin{align}
x(t+1)=Ax(t)+w(t),\, t\in\Z_{+},
\end{align}
where $x(t)\in\R^{n}$ is the state of the system and $w(t)$ is the
process noise. The initial state, $x(0)$, is assumed to be Gaussian
with zero mean and covariance matrix $\phi_{0}$, i.e., $x(0)\sim{\cal N}(0,\phi_{0})$.
There are $M$ different sensors attached to the process. At each
time step, we assume that only one sensor is available to take measurements.
The measurement of the $i$th sensor is given by: 
\begin{align}
y_{i}(t)=C_{i}x(t)+v_{i}(t),\, t\in\Z_{+},
\end{align}
where $y_{i}(t)\in\R^{p}$ and $v_{i}(t)\in\R^{p}$ are the measurement
output and measurement noise of the $i$th sensor at time $t$, respectively.
We assume that the process noise and all the measurement noises are
mutually independent Gaussian white noises given by: 
\begin{align*}
w(t)\sim{\cal N}(0,\Phi^{w}),\quad v_{i}(t)\sim{\cal N}(0,\Phi_{i}^{v}),
\end{align*}
all of which are also independent of the initial state $x(0)$.

Define $\lambda_{w}^{-}=\lambda_{\min}(\Phi^{w})$ and $\lambda_{v}^{-}=\min_{i\in\M}\{\lambda_{\min}(\Phi_{i}^{v})\}$.
Assume that $\lambda_{w}^{-}>0$ and $\lambda_{v}^{-}>0$. Let $\M:=\{1,\ldots,M\}$
be the set of sensor indices. For each $N\in\Z_{+}$, denote by $\M^{N}$
the set of all the sequences of sensor indices of length $N$. An
element $\sigma\in\M^{N}$ is called an $N$-horizon {\em sensor
schedule}. The set of all infinite-horizon sensor schedules is denoted
by $\M^{\infty}$. An infinite-horizon schedule $\sigma\in\M^{\infty}$
is called periodic with a period $l\in\Z_{+}$ if $\sigma(t)=\sigma(t+l)$
for all $t\in\Z_{+}$. Under a given sensor schedule $\sigma\in\M^{\infty}$,
the measurement sequence is determined by: 
\begin{align*}
y(t)=y_{\sigma(t)}(t)=C_{\sigma(t)}x(t)+v_{\sigma(t)}(t),\quad t\in\Z_{+}.
\end{align*}

For each $t_{1}\le t_{2}<\infty$, denote by $\hat{x}^{\sigma}(t_{2}|t_{1})$
the minimum mean-square error (MMSE) estimate of $x(t_{2})$ given
the measurements $\{y(0),\ldots,y(t_{1})\}$, the initial covariance
$\phi_{0}$ and the sensor schedule $\sigma\in\M^{\infty}$. Define
the prediction error $e^{\sigma}(t|t-1)$ by 
\begin{align*}
e^{\sigma}(t|t-1)=x(t)-A\hat{x}^{\sigma}(t-1|t-1),
\end{align*}
and let $\Sigma_{t}^{\sigma}(\phi_{0})$ be its covariance matrix.
When no ambiguity arises, we may drop the dependence on the initial
covariance matrix and simply write $\Sigma_{t}^{\sigma}$. By a standard
result of linear estimation theory, the prediction error covariance
can be updated recursively using the Riccati map:
\begin{eqnarray}
\Sigma_{t+1}^{\sigma} & = & \Phi^{w}+A\Sigma_{t}^{\sigma}A^{T}-A\Sigma_{t}^{\sigma}C_{\sigma(t)}^{T} \nonumber \\
 &  & \times\left(C_{\sigma(t)}\Sigma_{t}^{\sigma}C_{\sigma(t)}^{T}+\Phi_{\sigma(t)}^{v}\right)^{-1}C_{\sigma(t)}\Sigma_{t}^{\sigma}A^{T}.\label{eq:rr}
\end{eqnarray}

For any finite integer $N$, the performance of an $N$-horizon sensor
schedule $\sigma\in\M^{N}$ can be evaluated according to the total
estimation error defined by: 
\begin{align}
J_{N}(\sigma;\phi_{0})\triangleq\sum_{t=1}^{N}\tr(\Sigma_{t}^{\sigma}(\phi_{0})),\label{eq:JN}
\end{align}
or according to the average-per-stage estimation error defined by:
\begin{align}
\bar{J}_{N}(\sigma;\phi_{0})\triangleq\frac{1}{N}J_{N}(\sigma;\phi_{0}).\label{eq:JNbar}
\end{align}

When $N$ is finite, the two cost functions $J_{N}$ and $\bar{J}_{N}$
are equivalent in the sense that they produce the same set of optimal
solutions. However, the total cost $J_{N}(\sigma;\phi_{0})\to\infty$
as $N\to\infty$ for all $\sigma\in\M^{\infty}$ and $\phi_{0}\in{\cal A}$,
because the system is constantly perturbed by a nontrivial Gaussian
noise $w(t)$. Thus, the performance of an infinite-horizon sensor
schedule is usually measured by the limsup of the $N$-horizon average-per-stage
cost: 
\begin{align*}
\bar{J}_{\infty}(\sigma;\phi_{0})\triangleq\limsup_{N\to\infty}\bar{J}_{N}(\sigma;\phi_{0}).
\end{align*}

\begin{rem}
For notational convenience, the cost function is defined based on the prediction error covariance $\Sigma_{t}^{\sigma}\triangleq\text{cov}(e^\sigma(t|t-1))$ instead of the estimation error covariance (i.e. $\Sigma_{t|t}^{\sigma}\triangleq\text{cov}(e^\sigma(t|t))$). It can be shown that (pp. 80-81 of \cite{AbouKandil2003}) the iteration of $\Sigma_{t|t}^{\sigma}$ can be written in the same form as in (\ref{eq:rr}) with some modified system matrices ($\tilde{A}$, $\tilde{C}$), and process noise covariance $\tilde{\Phi}^w$. Therefore, all the results still hold if the estimation error covariance is used in the above cost function. 
\end{rem}

This cost function has been extensively used for studying various
infinite-horizon optimal control and estimation problems~\cite{BerDP01,NYACC09}.
However, this cost function depends only on the limiting behavior
of the schedule, which may lead to rather unexpected optimal solutions.
For example, one can manipulate a finite portion of an optimal schedule
to create an arbitrary transient behavior for the error trajectory
without affecting the optimality of the schedule. In some extreme
cases, the trajectory of the error covariance under an optimal schedule
may even grow unbounded. The following example illustrates such a
situation. 
\begin{example}
(Unbounded Optimal Schedule) Consider a simple 2-dimensional system
with $A=$diag~$\{\lambda,0\}$, $\lambda\in\Z_+$, $\lambda>1$, $C_{1}=[1,\ 0]$,
$C_{2}=[0,\ 1]$, $\Phi^{w}=\mbox{diag }\{c,c\}\succ0$, and $\Phi_{1}^{v}=\Phi_{2}^{v}=d>0.$ Note that the system
is detectable under sensor 1, but undetectable if using only sensor
2. One optimal schedule can be easily identified as using exclusively
sensor 1, which leads to a minimum cost of $\phi^{*}+c$, where $\phi^{*}=\frac{\sqrt{(d-c-\lambda^{2}d)^{2}+4cd}-(d-c-\lambda^{2}d)}{2}>c$ and diag\{$\phi^{*}$, $c$\} is the equilibrium point of the Riccati difference equation obtained from the covariance matrix iteration corresponding to sensor 1. Now
consider an infinite-horizon sensor schedule $\hat{\sigma}$ that
alternates between sensor 1 and sensor 2 on time intervals $\mathcal{I}_{k}^{1}$
and $\mathcal{I}_{k}^{2}$, respectively. Let the length of $\mathcal{I}_{k}^{1}$
be $k\lambda^{2k}$, and the length of $\mathcal{I}_{k}^{2}$ be $k$.
Define $t_{12}^{k}$ as the $k$th switching instants from sensor
2 to sensor 1 and $t_{21}^{k}$ the other way around. The switching
times can be determined recursively as: $t_{21}^{k}=t_{12}^{k}+k$
and $t_{12}^{k+1}=t_{21}^{k}+k\lambda^{2k}$, where $k=1,2,...$ and
$t_{12}^{1}=0$. Therefore, $\mathcal{I}_{k}^{2}=[t_{12}^{k},t_{21}^{k})$,
$\mathcal{I}_{k}^{1}=[t_{21}^{k},t_{12}^{k+1})$, and $\hat{\sigma}$
can be represented as 
\[
\hat{\sigma}(t)=\begin{cases}
2, & t\in[t_{12}^{k},t_{21}^{k})\\
1, & t\in[t_{21}^{k},t_{12}^{k+1}),
\end{cases}
\]
where $t\in\Z_{+}$, $t_{12}^{1}=0$ and $k=1,2,...$ It can be easily verified that the average-per-stage cost on $\cup_{i=1}^{k}\mathcal{I}_{k}$
with $\mathcal{I}_{k}=\mathcal{I}_{k}^{2}\cup\mathcal{I}_{k}^{1}$
goes to $\phi^{*}+c$ as $k\rightarrow\infty$, which is the same as the optimal
cost. However, the subsequence consisting of error covariances $\Sigma_{t_{21}^{k}-1}^{\sigma}$
diverges as $k\rightarrow\infty$. 
\end{example}
To exclude such abnormalities, we introduce the following feasible
set of sensor schedules which yield a bounded trajectory under a given
initial condition $\phi\in\mathcal{A}$: 
\[
\M_{\phi}^{\infty}=\{\sigma\in\M^{\infty}\!\!:\exists\beta<\infty,\text{ s.t. }\Sigma_{t}^{\sigma}(\phi)\!\preceq\!\beta I_{n},\forall t\in\Z_{+}\}.
\]

An infinite-horizon sensor schedule $\sigma$ is called {\em feasible}
for $\phi\in{\cal A}$ if $\sigma\in\M_{\phi}^{\infty}$. The following
assumption is adopted throughout this paper.

\begin{assumption} \label{ass:main} $\M_{\phi}^{\infty}\neq\emptyset$,
$\forall\phi\in{\cal A}$.

\end{assumption} 
\begin{rem}
The assumption requires that for any initial covariance, there always
exists an infinite-horizon schedule that can keep the estimation error
covariance bounded for all time. This is a reasonable assumption for
typical estimation applications. It can be guaranteed if, for example,
one of the subsystems is detectable. \end{rem}
\begin{problem}
\label{pro:infhss} For a given $\phi_{0}\in{\cal A}$, solve the
following problem 
\begin{equation}
\bar{V}^{*}(\phi_{0})\triangleq\inf_{\sigma\in\M_{\phi_{0}}^{\infty}}\limsup_{N\to\infty}\bar{J}_{N}(\sigma;\phi_{0})\label{eq:Vbarstar}
\end{equation}

\end{problem}
Assumption~\ref{ass:main} implies that $\bar{V}^{*}(\phi_{0})$
is finite for all $\phi_{0}\in{\cal A}$. The function $\bar{V}^{*}:{\cal A}\to\R_{+}$
defined implicitly by equation~(\ref{eq:Vbarstar}) is called the
optimal infinite-horizon cost function. For a general $\phi\in{\cal A}$,
a schedule that achieves the cost $V^{*}(\phi)$ will be referred
to as an {\em optimal schedule} for $\phi$.

\section{The Sequential Riccati Mapping and Its Stability \label{sec:srm}}

The Riccati recursion in~(\ref{eq:rr}) can be viewed as a mapping
from $\Sigma_{t}^{\sigma}\in\mathcal{A}$ to $\Sigma_{t+1}^{\sigma}\in\mathcal{A}$
depending on the sensor index chosen at time $t$. In general, for
each sensor $i\in\M$, we can define the {\em Riccati mapping}
as 
\begin{eqnarray*}
\rho_{i}(Q) & = & \Phi^{w}+AQA^{T}\\
 &  & -AQC_{i}^{T}\left(C_{i}QC_{i}^{T}+\Phi_{i}^{v}\right)^{-1}C_{i}QA^{T},\ Q\in{\cal A}.
\end{eqnarray*}

With this notation, for a generic initial covariance matrix $\phi\in{\cal A}$,
the covariance matrix $\Sigma_{t}^{\sigma}(\phi)$, defined in~(\ref{eq:rr}),
is the trajectory of the following matrix-valued time-varying nonlinear
system: 
\begin{align}
\Sigma_{t+1}^{\sigma}=\rho_{\sigma(t)}\left(\Sigma_{t}^{\sigma}\right),\,\text{ for }t\in\Z_{+},\text{ with }\Sigma_{0}^{\sigma}=\phi.\label{eq:Prm}
\end{align}

One can also view $\Sigma_{t}^{\sigma}(\cdot)$ as the composition
of a sequence of Riccati mappings, i.e. 
\begin{equation}
\Sigma_{t}^{\sigma}(\phi)=\rho_{\sigma(t-1)}\circ\rho_{\sigma(t-2)}\cdots\circ\rho_{\sigma(0)}(\phi),\quad t\in\Z_{+}.\label{eq:compoundR}
\end{equation}
We will also refer to $\Sigma_{t}^{\sigma}$ as the composite Riccati
map associated with $\sigma$.

To solve Problem~\ref{pro:infhss}, it is critical to understand
the dynamic behavior of the matrix-valued nonlinear system~(\ref{eq:Prm})
under different infinite-horizon schedules. Two well-known properties
of the Riccati mapping are useful for this purpose. 
\begin{lem}
\label{lem:rm} For any $i\in\M$, $Q_{1},Q_{2}\in{\cal A}$ and $c\in[0,1]$,
we have \end{lem}
\begin{enumerate}
\item $Q_{1}\preceq Q_{2}\Rightarrow\rho_{i}(Q_{1})\preceq\rho_{i}(Q_{2})$; 
\item $\rho_{i}(cQ_{1}+(1-c)Q_{2})\succeq c\rho_{i}(Q_{1})+(1-c)\rho_{i}(Q_{2})$. \end{enumerate}
\begin{rem}
The lemma indicates that the Riccati mapping is {\em monotone}
and {\em concave}. The monotonicity property is a well-known result
and its proof can be found in~\cite{KuSto86}. The concavity property
is an immediate consequence of Lemma 1-(e) in~\cite{BrTAC04}. 
\end{rem}
Based on these two properties, one can prove the following results. 
\begin{prop}
\label{prop:gbd}(Theorem 5 of~\cite{Journal_ss}) For any $\phi\in{\cal A}$,
$\epsilon\in\R^{+}$, $\sigma\in\M^{\infty}$, and $t\in\Z_{+}$,
we have $\Sigma_{t}^{\sigma}(\phi+\epsilon I_{n})\preceq\Sigma_{t}^{\sigma}(\phi)+g_{t}^{\sigma}(\phi)\cdot\epsilon$,
where $g_{t}^{\sigma}(\phi)$ is the directional derivative of the
$t$-step Riccati mapping $\Sigma_{t}^{\sigma}$ at $\phi$ along
direction $I_{n}$. Furthermore, if $\Sigma_{t}^{\sigma}(\phi)\preceq\beta I_{n}$
for all $t\in\Z_{+}$ and for some $\beta<\infty$, then $\tr(g_{t}^{\sigma}(\phi))\le n\beta\eta^{t}/\lambda_{w}^{-}$,
$\forall t\in\Z_{+}$, where 
\begin{align}
 & \eta=\frac{1}{1+\alpha\lambda_{w}^{-}}<1\quad\text{ and }\quad\alpha=\frac{\lambda_{w}^{-}}{\|A\|^{2}\beta^{2}+\lambda_{w}^{-}\beta}.\label{eq:alpha}
\end{align}

\end{prop}
The above theorem reveals an important property of system~(\ref{eq:Prm}),
namely, that boundedness of its trajectory implies an exponential
disturbance attenuation. This property plays a crucial role in the
derivation of the various properties of Problem~\ref{pro:infhss}
in Section~\ref{sec:pos}.

\section{Main Results \label{sec:pos}}

In this section, we will use the properties of the sequential Riccati
mapping derived in the last section to gain some insights into the
solution of Problem~\ref{pro:infhss}.

\subsection{Independence of the Initial Covariance\label{sec:ind}}

We first show that the feasible set is independent of the initial
covariance. 
\begin{lem}
\label{lem:fs} If $\sigma\in\M_{\phi_{1}}^{\infty}$ for some $\phi_{1}\in{\cal A}$,
then $\sigma\in\M_{\phi}^{\infty}$ for all $\phi\in{\cal A}$.\end{lem}
\begin{IEEEproof}
Fix arbitrary $\phi,\ \phi_{1}\in{\cal A}$, and $\sigma\in\M_{\phi_{1}}^{\infty}$.
Since $\phi\preceq\phi_{1}+\|\phi-\phi_{1}\|I_{n}$, we have by Proposition~\ref{prop:gbd},
\begin{align*}
\Sigma_{t}^{\sigma}(\phi) & \preceq\Sigma_{t}^{\sigma}(\phi_{1})+g_{t}^{\sigma}(\phi_{1})\cdot\|\phi-\phi_{1}\|.
\end{align*}

The first term on the right-hand side is bounded because $\sigma\in\M_{\phi_{1}}^{\infty}$,
while the second term is bounded due to Proposition~\ref{prop:gbd}.
Thus, $\sigma\in\M_{\phi}^{\infty}$. 
\end{IEEEproof}
Therefore, if an infinite-horizon schedule is feasible for some initial
covariance matrix, it is also feasible for all initial covariances.
This allows us to drop the dependence on the initial covariance and
simply define 
\begin{eqnarray*}
\M_{f}^{\infty} & = & \{\sigma\in\M^{\infty}:\forall\phi\in\mathcal{A},\ \exists\beta<\infty,\text{ s.t. }\\
 &  & \Sigma_{t}^{\sigma}(\phi)\preceq\beta I_{n},\forall t\in\Z_{+}\}.
\end{eqnarray*}

We next show that under a fixed schedule $\sigma\in\M_{f}^{\infty}$,
all the trajectories starting from different initial covariances will
eventually converge to the same trajectory. 
\begin{thm}
\label{thm:phiConv} For any feasible schedule $\sigma\in\M_{f}^{\infty}$,
we have that 
\[
\|\Sigma_{t}^{\sigma}(\phi_{1})-\Sigma_{t}^{\sigma}(\phi_{2})\|\to0\text{ exponentially as }t\to\infty,
\]
for all $\phi_{1},\phi_{2}\in{\cal A}$.\end{thm}
\begin{IEEEproof}
Fix arbitrary $\phi_{1}\in{\cal A}$ and $\phi_{2}\in{\cal A}$. Define
$\epsilon=\|\phi_{1}-\phi_{2}\|$. Without loss of generality, let
$\beta<\infty$ be the bound such that $\Sigma_{t}^{\sigma}(\phi_{i})\preceq\beta I_{n}$
for all $t\in\Z_{+}$ and $i=1,2$. By Proposition~\ref{prop:gbd},
we have 
\begin{eqnarray*}
\Sigma_{t}^{\sigma}(\phi_{2}) & \preceq & \Sigma_{t}^{\sigma}(\phi_{1}+\|\phi_{2}-\phi_{1}\|I_{n})\\
 & \preceq & \Sigma_{t}^{\sigma}(\phi_{1})+g_{t}^{\sigma}(\phi_{1})\cdot\epsilon\\
 & \preceq & \Sigma_{t}^{\sigma}(\phi_{1})+\left(\frac{n\beta\epsilon}{\lambda_{w}^{-}}\eta^{t}\right)\cdot I_{n}.
\end{eqnarray*}
Similarly, we can obtain $\Sigma_{t}^{\sigma}(\phi_{1})\preceq\Sigma_{t}^{\sigma}(\phi_{2})+\left(\frac{n\beta\epsilon}{\lambda_{w}^{-}}\eta^{t}\right)\cdot I_{n},$
for all $t\in\Z_{+}$. The result follows directly from the above
inequalities as $t\to\infty$. 
\end{IEEEproof}
An immediate consequence of the above theorem is that the infinite-horizon
average-per-stage cost of any feasible schedule is independent of
the initial covariance matrix. 
\begin{cor}
\label{cor:Jinfphi} For any $\sigma\in\M_{f}^{\infty}$, $\bar{J}_{\infty}(\sigma;\phi_{1})=\bar{J}_{\infty}(\sigma;\phi_{2})$
for all $\phi_{1},\phi_{2}\in{\cal A}$.\end{cor}
\begin{IEEEproof}
By Theorem~\ref{thm:phiConv}, $\Sigma_{t}^{\sigma}(\phi_{1})\to\Sigma_{t}^{\sigma}(\phi_{2})$
as $t\to\infty$. Thus, the two average-per-stage cost sequences $\{\tfrac{1}{N}\sum_{t=1}^{N}\Sigma_{t}^{\sigma}(\phi_{i})\}_{N}$,
$N\in\Z_{+}$, $i=1,2${\small{,}} must have the same limsup. 
\end{IEEEproof}
By the above corollary, it is easy to see that if a feasible schedule
$\sigma$ is optimal for some initial covariance $\phi_{1}$, then
it must also be optimal for any other initial covariance $\phi_{2}$.
In addition, the optimal infinite-horizon average-per-stage costs
corresponding to these two initial covariances must also be the same. 
\begin{cor}
\label{cor:vstarphi} For any $\phi_{1},\phi_{2}\in{\cal A}$, if
$\sigma^{*}$ is optimal for $\phi_{1}$, then it must also be optimal
for $\phi_{2}$; and in addition, $\bar{V}^{*}(\phi_{1})=\bar{V}^{*}(\phi_{2})$. 
\end{cor}
Therefore, to solve Problem~\ref{pro:infhss}, we can start from
any initial covariance matrix at our convenience. The obtained optimal
solution would also be optimal for all the other initial covariances.

\subsection{Properties of the Accumulation Set}

For any $\sigma\in\M_{f}^{\infty}$, let ${\cal L}^{\sigma}$ be the
{\em accumulation set} of the closed-loop trajectory of the nonlinear
system~(\ref{eq:Prm}) under schedule $\sigma$. In other words,
the set ${\cal L}^{\sigma}$ contains all the points whose arbitrary
neighborhoods will be visited infinitely often by the trajectory $\{\Sigma_{t}^{\sigma}(\phi)\}_{t\in\Z_{+}}$
for some initial condition $\phi\in{\cal A}$. Under Assumption \ref{ass:main},
the sequence $\{\Sigma_{t}^{\sigma}(\phi)\}_{t\in\Z_{+}}$ is bounded
if $\sigma$ is feasible. Therefore, there exists a convergent subsequence
and ${\cal L}^{\sigma}$ is not empty. Moreover, ${\cal L}^{\sigma}$
is closed since the subsequential limits of a bounded sequence in
a metric space $X$ form a closed subset of $X$. It follows that
${\cal L}^{\sigma}$ is bounded and closed in ${\cal A}$, and is
thus compact.

According to Theorem~\ref{thm:phiConv}, a trajectory $\{\Sigma_{t}^{\sigma}(\phi)\}_{t\in\Z_{+}}$
under schedule $\sigma$ starting from any initial covariance $\phi\in{\cal A}$
has the same accumulation set ${\cal L}^{\sigma}$. This implies that
${\cal L}^{\sigma}$ is globally attractive, i.e. $\underset{t\to\infty}{\lim}d(\Sigma_{t}^{\sigma}(\phi),{\cal L}^{\sigma})=0,\forall\phi\in{\cal A}$
, where $d(\phi,{\cal L}^{\sigma})=\underset{z\in{\cal L}^{\sigma}}{\inf}\left\Vert \phi-z\right\Vert $
represents the distance from the point $\phi$ to the set ${\cal L}^{\sigma}$.

We summarize the above results in the following proposition: 
\begin{prop}
The accumulation set of any feasible schedule is nonempty, compact
and globally attractive. 
\end{prop}

\subsection{Universal Approximation Property of Periodic Schedules \label{sub:UAP}}

The goal of this subsection is to show that the optimal infinite-horizon
cost can be approximated with an arbitrary accuracy by periodic schedules.
Actually a more general result is proved for approximating infinite-horizon
costs of any feasible schedule. First, we derive the following result
which will facilitate our main proof. 
\begin{lem}
(Uniform Bound) \label{lem:unifBd} Given $\sigma\in\M_{f}^{\infty}$,
for any bounded set $E\subset{\cal A}$, there exists finite constants
$\beta_{E}$, $\alpha_{E}$, and $\eta_{E}\in(0,1)$, such that $\Sigma_{t}^{\sigma}(\phi)\preceq\beta_{E}I_{n}$
and $\tr(g_{t}^{\sigma}(\phi))\le\alpha_{E}\eta_{E}^{t}$, for all
$t\in\Z_{+}$ and $\phi\in E$.\end{lem}
\begin{IEEEproof}
Fix an arbitrary $\phi_{1}\in E$. Define the covariance trajectory
under $\sigma$ with initial covariance $\phi_{1}$ as $\psi_{t}=\Sigma_{t}^{\sigma}(\phi_{1})$,
$t\in\Z_{+}$. Since $\sigma$ is feasible, there must exist a finite
constant $\beta_{1}$ such that $\psi_{t}\le\beta_{1}I_{n}$ for all
$t\in\Z_{+}$. By Proposition~\ref{prop:gbd}, there exist constants
$\alpha_{1}<\infty$ and $\eta_{1}\in(0,1)$ such that $\tr(g_{t}^{\sigma}(\phi_{1}))\le\alpha_{1}\eta_{1}^{t}$,
for all $t\in\Z_{+}$. It follows that 
\begin{align*}
\Sigma_{t}^{\sigma}(\phi) & \preceq\Sigma_{t}^{\sigma}(\phi_{1}+\|\phi-\phi_{1}\|I_{n}) \\
& \preceq\Sigma_{t}^{\sigma}(\phi_{1})+g_{t}^{\sigma}(\phi_{1})\|\phi-\phi_{1}\|\\
& \preceq\psi_{t}+\alpha_{1}\eta_{1}^{t}(\kappa_{E}+\beta_{1})I_{n} \\
& \preceq\left[\beta_{1}+\alpha_{1}\eta_{1}^{t}(\kappa_{E}+\beta_{1})\right]I_{n} \\
& \triangleq\beta_{E}I_{n},
\end{align*}
for all $\phi\in E$, where $\kappa_{E}\triangleq\sup_{\phi\in E}\|\phi\|$.
This implies the existence of the desired constant $\beta_{E}$, which is common for all trajectories starting from $E$. This in turn guarantees the existence of the desired constants $\alpha_{E}$
and $\eta_{E}$ according to Proposition~\ref{prop:gbd}. 
\end{IEEEproof}
The above lemma indicates that the covariance trajectories starting
from any initial covariance in a bounded set $E$ are bounded uniformly
by $\beta_{E}I_{n}$. The bound $\beta_{E}$ depends only on the underlying
set $E$ instead of the particular value of the initial covariance.

The following theorem presents the main contribution of this paper. 
\begin{thm}
(Universal Approximation)\label{thm:pss} For any feasible schedule
$\sigma\in\M_{f}^{\infty}$ and any $\delta>0$, there exists a periodic
schedule $\tilde{\sigma}$ with a finite period $N\in\Z_{+}$, such
that the infinite-horizon cost of $\sigma$ is approximated by $\tilde{\sigma}$
with the error bound $\left|\bar{J}_{\infty}(\tilde{\sigma})-\bar{J}_{\infty}(\sigma)\right|<\delta.$
\end{thm}
%The proof proceeds in four steps: 1) define the domain $E$ of consideration;
%2) find a finite schedule $\sigma_{l}$ which uniformly steers all
%the trajectories starting from $E$ into a small neighborhood of an
%accumulation point; 3) find finite schedules $\sigma_{N_{k}}$ whose
%average-per-stage costs converge to the infinite-horizon cost of $\sigma$
%uniformly for all initial condition in $E$; 4) construct the periodic
%schedule $\tilde{\sigma}$ out of $\sigma_{N_{k}} and \sigma_{l}$, and
%prove that its infinite cost satisfies the error bound $\delta$ for
%all large enough $k$.

\begin{IEEEproof}
Pick an arbitrary feasible schedule $\sigma\in\M_{f}^{\infty}$ and
an accumulation point $\hat{\phi}\in\mathcal{L^{\sigma}}$. Suppose
that $\Sigma_{t}^{\sigma}(\hat{\phi})<\beta I_{n}$. By Proposition~\ref{prop:gbd}
and Lemma \ref{lem:unifBd}, we have 
\[
\Sigma_{t}^{\sigma}(\phi)\preceq\Sigma_{t}^{\sigma}(\hat{\phi})+r\alpha_{r}\eta_{r}^{t}I_{n}\preceq(\beta+r\alpha_{r}\eta_{r}^{t})I_{n},\quad\forall\phi\in\mathcal{B}(\hat{\phi};r)
\]
where $\alpha_{r}>0$ and $0<\eta_{r}<1$ are constants depending
on $r$. Denote $\hat{\beta}=\beta+r\alpha_{r}$. It is clear that
$\Sigma_{t}^{\sigma}(\phi)$ is bounded by $\hat{\beta}I_{n}$, $\forall\phi\in\mathcal{B}(\hat{\phi};r)$.
Define $E=\{\phi:\ \phi\preceq\hat{\beta}I_{n}\}$. Clearly $\mathcal{L^{\sigma}}\subset E$,
and $\mathcal{B}(\hat{\phi};r)\subset E$. These sets are illustrated
in Fig.~\ref{fig:DC}.

The rest of the proof consists of three major steps. (i) Firstly,
we show that there exists a common $l$-horizon schedule $\sigma_{l}$
that can drive the covariance trajectory to $\mathcal{B}(\hat{\phi};r)$
at the end of the $l$ horizon for any initial covariance in $E$;
(ii) Secondly, we show that there exists a subschedule $\sigma_{N_{k}}$
whose average-per-stage cost converges to the infinite-horizon cost
of $\sigma$ uniformly for all initial condition in $E$; (iii) Lastly,
we will construct a periodic schedule $\tilde{\sigma}$ based on $\sigma_{l}$
and $\sigma_{N_{k}}$, which satisfies the desired error bound $\delta$
for all large enough $k$. 
\begin{figure}[t]
\centering{}\includegraphics[clip,width=0.35\textwidth]{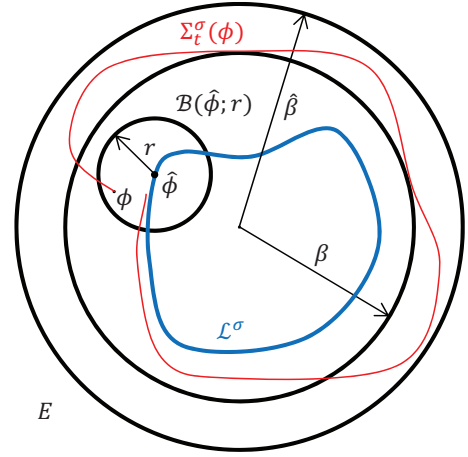}\caption{\label{fig:DC}Domain of consideration: bounded set $E$ which contains
all the trajectories starting from $\mathcal{B}(\hat{\phi};r)$}
\end{figure}

\textbf{Step (i): }By Proposition~\ref{prop:gbd} and Lemma \ref{lem:unifBd},
$\forall\phi\in E$, {\small{
\begin{equation}
\left\Vert \Sigma_{t}^{\sigma}(\phi)-\Sigma_{t}^{\sigma}(\hat{\phi})\right\Vert\leq\hat{\beta}\alpha_{E}\eta_{E}^{t},\label{eq:init}
\end{equation}
}}where $\alpha_{E}>0$ and $0<\eta_{E}<1$ are constants associated
with $E$. Therefore $\exists l_{0}>0$ such that\textbf{ }\textbf{\small{
\begin{equation}
\left\Vert \Sigma_{t}^{\sigma}(\phi)-\Sigma_{t}^{\sigma}(\hat{\phi})\right\Vert\leq\frac{r}{2},\ \forall t>l_{0}.\label{eq:close}
\end{equation}
}}{\small \par}

Since $\hat{\phi}\in\mathcal{L^{\sigma}}$ and $\mathcal{L^{\sigma}}$
is attractive, $\exists l\geq l_{0}$ such that {\small{
\begin{equation}
\left\Vert \Sigma_{l}^{\sigma}(\hat{\phi})-\hat{\phi}\right\Vert \leq\frac{r}{2}.\label{eq:attractive}
\end{equation}
}}{\small \par}

By~(\ref{eq:close}) and~(\ref{eq:attractive}), we have $\forall\phi\in E${\small{
\begin{equation}
\left\Vert \Sigma_{l}^{\sigma}(\phi)-\hat{\phi}\right\Vert \leq\left\Vert \Sigma_{l}^{\sigma}(\phi)-\Sigma_{l}^{\sigma}(\hat{\phi})\right\Vert +\left\Vert \Sigma_{l}^{\sigma}(\hat{\phi})-\hat{\phi}\right\Vert \leq r.\label{eq:newcontraction}
\end{equation}
}}{\small \par}

Denote the first $l$ steps of $\sigma$ by $\sigma_{l}$. Equation~(\ref{eq:newcontraction})
shows that under $\sigma_{l}$, the final covariance $\Sigma_{l}^{\sigma_{l}}(\phi)\in{\cal B}(\hat{\phi};r)$,
$\forall\phi\in E$.

\textbf{Step (ii): }Now we construct another finite length sub-schedule
$\sigma_{N_{k}}$, under which the performance obtained can be arbitrarily
close to $\bar{J}_{\infty}(\sigma)$ when $k$ is large enough.

Suppose $\bar{J}_{N}(\sigma;\phi_{0})$ is the average-per-stage cost
of the first $N$ steps of $\sigma$ for some initial condition $\phi_{0}\in E$.
For brevity, define $b_{N}=\bar{J}_{N}(\sigma;\phi_{0})$. Since $\sigma$
is feasible, $b_{N}$ is bounded. Therefore, there exists a subsequence
$\left\{ b_{N_{k}}\right\} _{k}$, such that $\bar{J}_{\infty}(\sigma;\phi_{0})=\underset{{\scriptstyle k\rightarrow\infty}}{\lim}b_{N_{k}}.$
It follows $\forall\delta>0$, $\exists K_{1}\in\mathbb{Z}_{+}$,
such that{\small{
\begin{equation}
\left|\bar{J}_{\infty}(\sigma;\phi_{0})-b_{N_{k}}\right|<\frac{\delta}{3},\ \ \forall k>K_{1}.\label{eq:Finite}
\end{equation}
}}{\small \par}

Let $\sigma_{N_{k}}$ be the first $N_{k}$ steps of $\sigma$. Note
that $\left\{ b_{N_{k}}\right\} _{k}$ and the associated sub-schedule
$\{\sigma_{N_{k}}\}_{k}$ are constructed based on the initial condition
$\phi_{0}$. We need to be shown that the cost convergence is uniform
with respect to all $\phi\in E$. To this end, consider arbitrary
$\phi_{1},\phi_{2}\in E$. We know that $\left\Vert \phi_{1}-\phi_{2}\right\Vert \leq\hat{\beta}$.
It follows from Proposition~\ref{prop:gbd} and Lemma \ref{lem:unifBd}
that {\small{
\begin{eqnarray*}
\Delta & \triangleq & \left|\bar{J}_{N_{k}}(\sigma_{N_{k}};\phi_{1})-\bar{J}_{N_{k}}(\sigma_{N_{k}};\phi_{2})\right|\\
 & \leq & \frac{1}{N_{k}}\overset{{\scriptstyle N_{k}}}{\underset{{\scriptstyle t=1}}{\Sigma}}\left|\mbox{tr}(\Sigma_{t}^{\sigma}(\phi_{1}))-\mbox{tr}(\Sigma_{t}^{\sigma}(\phi_{2}))\right|\\
 & \leq & \frac{1}{N_{k}}\overset{{\scriptstyle N_{k}}}{\underset{{\scriptstyle t=1}}{\Sigma}}n\hat{\beta}\alpha_{E}\eta_{E}^{t}\leq\frac{n\hat{\beta}\alpha_{E}\eta_{E}}{N_{k}(1-\eta_{E})}.
\end{eqnarray*}
}}{\small \par}

Note that $n$, $\hat{\beta}$, $\alpha_{E}$, $\eta_{E}$ are all
constants. Hence, $\exists K_{2}\in\mathbb{Z}_{+}$, such that $\Delta<\frac{\delta}{3}$,
$\forall k>K_{2}$. Choose $K=\max\{K_{1},K_{2}\}$, we have {\small{
\begin{eqnarray*}
\left|\bar{J}_{\infty}(\sigma;\phi_{0})-\bar{J}_{N_{k}}(\sigma_{N_{k}};\phi)\right| & \leq & \left|\bar{J}_{\infty}(\sigma;\phi_{0})-\bar{J}_{N_{k}}(\sigma_{N_{k}};\phi_{0})\right|\\
 &  & +\left|\bar{J}_{N_{k}}(\sigma_{N_{k}};\phi_{0})-\bar{J}_{N_{k}}(\sigma_{N_{k}};\phi)\right|\\
 & < & \frac{2\delta}{3},\ \ \forall\phi\in E,\ \ \forall k>K
\end{eqnarray*}
}}{\small \par}

\textbf{Step (iii): }Now construct a periodic schedule $\tilde{\sigma}=\{\sigma_{N},\sigma_{N},\cdots\cdots\}$,
where $\sigma_{N}=\{\sigma_{N_{k}},\sigma_{l}\}$, $N=N_{k}+l$. Recall
that $N_{k}$ is constructed in {\em Step (ii)} for some large
$k$ to be determined later; and $l\geq l_{0}$ is from {\em Step
(i)} so that equations (\ref{eq:attractive}) and (\ref{eq:newcontraction})
hold. Note that $\Sigma_{N_{k}}^{\sigma_{N_{k}}}(\phi)\in E$, $\forall\phi\in\mathcal{B}(\hat{\phi};r)$.
It follows from~(\ref{eq:newcontraction}) that $\Sigma_{N}^{\sigma_{N}}(\phi)\in\mathcal{B}(\hat{\phi};r)$,
$\forall\phi\in\mathcal{B}(\hat{\phi};r)$. Therefore $\Sigma_{N}^{\sigma_{N}}(\phi)$
is an invariant mapping with respect to $\mathcal{B}(\hat{\phi};r)$.
Further note that~(\ref{eq:init}) and~(\ref{eq:close}) implies
$\forall\phi_{1},\phi_{2}\in\mathcal{B}(\hat{\phi};r)$, 
\begin{eqnarray*}
\left\Vert \Sigma_{t}^{\sigma}(\phi_{1})-\Sigma_{t}^{\sigma}(\phi_{2})\right\Vert  & \leq & \alpha_{E}\eta_{E}^{t}\left\Vert \phi_{1}-\phi_{2}\right\Vert \\
 & \leq & \frac{r}{2\hat{\beta}}\left\Vert \phi_{1}-\phi_{2}\right\Vert ,\ \ \forall t>l_{0}
\end{eqnarray*}
Recall that $\mathcal{B}(\hat{\phi};r)\subset E$ as discussed in
Step 1, and therefore $r<2\hat{\beta}$. It follows that $\Sigma_{N}^{\sigma_{N}}$
is a contraction mapping on $\mathcal{B}(\hat{\phi};r)$. Let $P$
be the unique fixed point of $\Sigma_{N}^{\sigma_{N}}$ on $\mathcal{B}(\hat{\phi};r)$.
Since $P\in\mathcal{L}^{\tilde{\sigma}}$ and $\mathcal{L}^{\tilde{\sigma}}$
is an $N$-cycle %
\footnote{By $N$-cycle, we mean a sequence $(\phi_{0},\phi_{1},\cdots,\phi_{N-1})$
where $\rho_{\tilde{\sigma}(0)}(\phi_{0})=\phi_{1}$, $\rho_{\tilde{\sigma}(1)}(\phi_{1})=\phi_{2}$,$\cdots$,
and $\rho_{\tilde{\sigma}(N-1)}(\phi_{N-1})=\phi_{0}$.%
}, the performance of $\tilde{\sigma}$ can be obtained as 
\begin{eqnarray*}
\bar{J}_{\infty}(\tilde{\sigma};P) & = & \frac{N_{k}\bar{J}_{N_{k}}(\sigma_{N_{k}};P)+l\bar{J}_{l}(\sigma_{l};\phi_{P})}{N_{k}+l},
\end{eqnarray*}
where $\phi_{P}$ denotes the point $\Sigma_{N_{k}}^{\sigma_{N_{k}}}(P)$.

Since $l$ is independent of $N_{k}$ and $\bar{J}_{l}(\sigma_{l};\phi_{P})$
is bounded, $\bar{J}_{N_{k}}(\sigma_{N_{k}};P)\to\bar{J}_{\infty}(\tilde{\sigma};P)$
as $k\rightarrow\infty$. That is $\forall\delta>0,$ $\exists K_{3}\in\mathbb{Z}_{+}$,
such that $\left|\bar{J}_{\infty}(\tilde{\sigma};P)-\bar{J}_{N_{k}}(\sigma_{N_{k}};P)\right|<\frac{\delta}{3}$,
$\forall k>K_{3}$. Let $\bar{K}=\max\{K_{3},K\}$, and it follows
$\left|\bar{J}_{\infty}(\tilde{\sigma};P)-\bar{J}_{\infty}(\sigma;P)\right|<\delta$
when $k>\bar{K}$ and the length of the period of $\tilde{\sigma}$
is $N=N_{\bar{K}}+l$. By Theorem~\ref{thm:phiConv}, we know that
the infinite-horizon cost is independent of the initial condition,
and thus the desired result follows. \end{IEEEproof}
\begin{rem}
\label{rmk:subseq}The proof of Theorem~\ref{thm:pss} can also be
applied to prove that the $\liminf$ (or any subsequential limits)
of the sequence of average-per-stage costs can be approximated by
periodic schedules arbitrarily closely by choosing $\sigma_{N_{k}}$
corresponding to the $\liminf$ subsequences (or any convergent subsequences).
Then the following corollary follows immediately from Theorem~\ref{thm:pss}:\end{rem}
\begin{cor}
\label{cor:Conv}Suppose $\sigma^{*}$ is an optimal sensor schedule
to Problem~\ref{pro:infhss}, and $\left\{ b_{N}^{*}\right\} $ is
the corresponding sequence of the $N$-horizon average-per-stage costs.
Then $\left\{ b_{N}^{*}\right\} $ converges and the optimal cost
is $\bar{V}^{*}=$$\underset{N\rightarrow\infty}{\lim}b_{N}^{*}$.\end{cor}
\begin{IEEEproof}
Suppose $\left\{ b_{N}^{*}\right\} $ does not converge. Let $\varepsilon=\underset{N\rightarrow\infty}{\limsup}b_{N}^{*}-\underset{N\rightarrow\infty}{\liminf}b_{N}^{*}>0$.
By Theorem~\ref{thm:pss} and Remark~\ref{rmk:subseq}, there exists
a periodic schedule $\tilde{\sigma}$ with finite period such that
{\small{$\left|\bar{J}_{\infty}(\tilde{\sigma})-\underset{N\rightarrow\infty}{\liminf}b_{N}^{*}\right|<\frac{\varepsilon}{2}$}}.
It follows $\bar{J}_{\infty}(\tilde{\sigma})<\underset{N\rightarrow\infty}{\liminf}b_{N}^{*}+\frac{\varepsilon}{2}<\underset{N\rightarrow\infty}{\limsup}b_{N}^{*}$.
Thus $\tilde{\sigma}$ has a smaller infinite cost than the optimal
schedule $\sigma^{*}$, which is a contradiction. 
\end{IEEEproof}
In addition, the proof of Theorem~\ref{thm:pss} also implies the
stability of the covariance trajectory under a feasible periodic sensor
schedule. 
\begin{cor}
For any feasible periodic schedule $\tilde{\sigma}$, the discrete
nonlinear system $\phi_{k+1}=\Sigma_{N}^{\tilde{\sigma}}(\phi_{k})$,
$k\in\Z_{+}$, $\forall\phi_{0}\in\mathcal{A}$ is globally asymptotically
and locally exponentially stable.\end{cor}
\begin{IEEEproof}
From the proof of Theorem~\ref{thm:pss}, the fixed point $P$ of
the contraction mapping $\Sigma_{N}^{\tilde{\sigma}}$ on $\mathcal{B}(\hat{\phi};r)$
is also the equilibrium of the system $\phi_{k+1}=\Sigma_{N}^{\tilde{\sigma}}(\phi_{k})$.
The result follows by futher noting that $P\in{\cal L}^{\tilde{\sigma}}$
is globally attractive. 
\end{IEEEproof}

\section{Discussion and Conclusions \label{sec:con}}

Under a mild feasibility assumption, we have proven
that both the optimal infinite-horizon cost and the corresponding
optimal schedule are independent of the initial error covariance.
Furthermore, we have proven that the accumulation set of the composite Riccati
mapping under a feasible schedule is nonempty, compact, and globally
attractive. The most important result is the universal approximation
theorem (Theorem~\ref{thm:pss}) which states that the performance
of any feasible schedule can be approximated arbitrarily closely
by a periodic schedule with a finite period. Interestingly, this result
leads to the conclusion that the average-per-stage cost of an optimal
schedule must converge (Corollary~\ref{cor:Conv}).

These theoretical results provide us valuable insights into the infinite-horizon
sensor scheduling problem. Theorem~\ref{thm:pss} motivates us to
focus on the periodic schedules in solving Problem~\ref{pro:infhss}. Corollary~\ref{cor:Conv} can be used to simplify the objective function
(cost function) in optimization-based approaches for finding the optimal/suboptimal
periodic solutions. It is worth mentioning that there are many other ways to quantify the performance of an infinite-horizon sensor schedule. For example, one can also use a discounted cost function $\sum_{t=1}^\infty c^t\text{tr}(\Sigma^\sigma_t(\phi))$  with a discount factor $c\in (0,1)$ or simply use the limsup of the terminal cost $\underset{N\rightarrow\infty}{\limsup}\ tr(\Sigma_{N}^{\sigma})$. In general, our results do not directly apply to these cost metrics. Our future research will focus on extending the results to other cost metrics, and on developing efficient infinite-horizon sensor scheduling algorithms with guaranteed suboptimal performance.

\bibliographystyle{IEEEtranS}
\bibliography{ss}

% Generated by IEEEtranS.bst, version: 1.13 (2008/09/30)
\begin{thebibliography}{10}
\providecommand{\url}[1]{#1}
\csname url@samestyle\endcsname
\providecommand{\newblock}{\relax}
\providecommand{\bibinfo}[2]{#2}
\providecommand{\BIBentrySTDinterwordspacing}{\spaceskip=0pt\relax}
\providecommand{\BIBentryALTinterwordstretchfactor}{4}
\providecommand{\BIBentryALTinterwordspacing}{\spaceskip=\fontdimen2\font plus
\BIBentryALTinterwordstretchfactor\fontdimen3\font minus
  \fontdimen4\font\relax}
\providecommand{\BIBforeignlanguage}[2]{{%
\expandafter\ifx\csname l@#1\endcsname\relax
\typeout{** WARNING: IEEEtranS.bst: No hyphenation pattern has been}%
\typeout{** loaded for the language `#1'. Using the pattern for}%
\typeout{** the default language instead.}%
\else
\language=\csname l@#1\endcsname
\fi
#2}}
\providecommand{\BIBdecl}{\relax}
\BIBdecl

\bibitem{AbouKandil2003}
H.~Abou-Kandil, G.~Freiling, V.~Ionescu, and G.~Jank, \emph{Matrix Riccati
  Equations in Control and Systems Theory}.\hskip 1em plus 0.5em minus
  0.4em\relax Birkh{\"{a}}user Verlag, 2003.

\bibitem{BeAUTO05}
S.~C. Bengea and R.~A. DeCarlo, ``Optimal control of switching systems,''
  \emph{Automatica}, vol.~41, no.~1, pp. 11--27, 2005.

\bibitem{BerDP01}
D.~Bertsekas, \emph{Dynamic Programming and Optimal Control}, 2nd~ed.\hskip 1em
  plus 0.5em minus 0.4em\relax Athena Scientific, 2001, vol.~2.

\bibitem{FEAUTO08}
Z.~G. Feng, K.~L. Teo, and V.~Rehbock, ``Hybrid method for a general optimal
  sensor scheduling problem in discrete time,'' \emph{Automatica}, vol.~44,
  no.~5, pp. 1295--1303, 2008.

\bibitem{GupAuto06}
V.~Gupta, T.~H. Chung, B.~Hassibi, and R.~M. Murray, ``On a stochastic sensor
  selection algorithm with applications in sensor scheduling and sensor
  coverage.'' \emph{Automatica}, vol.~42, no.~2, pp. 251--260, 2006.

\bibitem{HoCDC07}
P.~Hovareshti, V.~Gupta, and J.~S. Baras, ``Sensor scheduling using smart
  sensors,'' in \emph{IEEE Conference on Decision and Control}, New Orleans,
  LA, Dec. 2007.

\bibitem{JOTSP09}
S.~Joshi and S.~Boyd, ``Sensor selection via convex optimization,'' \emph{IEEE
  Transactions on Signal Processing}, vol.~57, no.~2, pp. 451--462, 2009.

\bibitem{KuSto86}
P.~R. Kumar and P.~Varaiya, \emph{Stochastic systems: estimation,
  identification and adaptive control}.\hskip 1em plus 0.5em minus 0.4em\relax
  Upper Saddle River, NJ, USA: Prentice-Hall, Inc., 1986.

\bibitem{Lemmon2010}
M.~Lemmon, ``Event-triggered feedback in control, estimation, and
  optimization,'' in \emph{Lecture Notes in Control and Information Sciences},
  A.~Bemporad, M.~Heemels, and M.~Johansson, Eds.\hskip 1em plus 0.5em minus
  0.4em\relax Springer London, 2010, vol. 406, pp. 293--358.

\bibitem{LiDSP09}
Y.~Li, L.~Krakow, E.~Chong, and K.~Groom, ``Approximate stochastic dynamic
  programming for sensor scheduling to track multiple targets,'' \emph{Digital
  Signal Processing}, vol.~19, no.~6, pp. 978--989, Dec. 2009.

\bibitem{LOACC99}
A.~Logothetis and A.~Isaksson, ``On sensor scheduling via information theoretic
  criteria,'' in \emph{Proceedings of the American Control Conference}, San
  Diego, CA, Jun. 1999.

\bibitem{MeiTAC67}
L.~Meier, J.~Peschon, and R.~Dressler, ``Optimal control of measurement
  subsystems,'' \emph{IEEE Transactions on Automatic Control}, vol.~12, no.~5,
  pp. 528--536, 1967.

\bibitem{MOTR06}
A.~I. Mourikis and S.~I. Roumeliotis, ``Optimal sensor scheduling for
  resource-constrained localization of mobile robot formations,'' \emph{IEEE
  Transactions on Robotics}, vol.~22, no.~5, pp. 917 -- 931, 2006.

\bibitem{NYACC09}
J.~L. Ny, E.~Feron, and M.~A. Dahleh, ``Scheduling {Kalman} filters in
  continuous time,'' in \emph{Proceedings of the American Control Conference},
  St Louis, MO, Jun. 2009.

\bibitem{OSTAES}
Y.~Oshman, ``Optimal sensor selection strategy for discrete-time state
  estimators,'' \emph{IEEE Transactions on Aerospace and Electronic Systems},
  vol.~30, no.~2, pp. 307--314, 1994.

\bibitem{Shi2011}
L.~Shi, Johansson, K.~Henrik, and L.~Qiu, ``Time and event-based sensor
  scheduling for networks with limited communication resources,'' in
  \emph{Proceedings of the 18th IFAC World Congress}, 2011, pp.
  13\,263--13\,268.

\bibitem{Sijs2009}
J.~Sijs and M.~Lazar, ``On event based state estimation,'' in \emph{Lecture
  Notes in Computer Science}, R.~Majumdar and P.~Tabuada, Eds.\hskip 1em plus
  0.5em minus 0.4em\relax Springer Berlin Heidelberg, 2009, vol. 5469, pp.
  336--350.

\bibitem{BrTAC04}
B.~Sinopoli, L.~Schenato, M.~Franceschetti, K.~Poolla, M.~I. Jordan, and S.~S.
  Sastry, ``{Kalman} filtering with intermittent observations,'' \emph{IEEE
  Transactions on Automatic Control}, vol.~49, no.~9, pp. 1453--1464, Sep.
  2004.

\bibitem{ACC10_ss}
M.~P. Vitus, W.~Zhang, A.~Abate, J.~Hu, and C.~J. Tomlin, ``On efficient sensor
  scheduling for linear dynamical systems,'' in \emph{Proceedings of the
  American Control Conference}, Baltimore, MD, Jul. 2010.

\bibitem{CDC10_Mike}
------, ``On sensor scheduling of linear dynamical systems with error bounds,''
  in \emph{IEEE Conference on Decision and Control}, Atlanta, GA, Dec. 2010.

\bibitem{Journal_ss}
------, ``On efficient sensor scheduling for linear dynamical systems,''
  \emph{Automatica}, vol.~48, no.~10, pp. 2482--2493, Oct. 2012.

\end{thebibliography}

\end{document}